\documentclass{amsart}
\usepackage{egen, egenengelsk}
\usepackage[twoside,left=3.66cm,right=3.66cm,bottom=3cm,top=2.66cm]{geometry}
\usepackage{hyperref}
\usepackage{tikz-cd}

\DeclareMathOperator{\GL}{GL}
\DeclareMathOperator{\Span}{Span}
\DeclareMathOperator{\td}{td}
\DeclareMathOperator{\Sing}{Sing}

\DeclareMathOperator{\PGL}{PGL}
\DeclareMathOperator{\SL}{SL}

\newcommand{\gl}{\mathfrak{gl}}
\newcommand{\op}{\text{op}}
\newcommand{\cald}{Căldăraru}
\newcommand{\szendroi}{Szendrői}
\newcommand{\ttop}{\text{top}}
\newcommand{\odd}{\text{odd}}
\renewcommand{\SS}{\ensuremath{\mathbb{S}}}

\newtheorem*{notbirational}{Theorem \ref{thm:XAndYNotBirational}}

\title[Counterexample to the birational Torelli problem for Calabi--Yau 3-folds]{A counterexample to the birational Torelli problem for Calabi--Yau threefolds}
\author{John Christian Ottem}
\address[Ottem]{Department of Mathematics \\
University of Oslo \\
PO Box 1053, Blindern \\
0316 Oslo \\
Norway}
\email{johnco@math.uio.no}

\author{Jørgen Vold Rennemo}
\address[Rennemo]{All Souls College \\
Oxford \\
OX1 4AL \\
UK}
\curraddr{Department of Mathematics \\
University of Oslo \\
PO Box 1053, Blindern \\
0316 Oslo \\
Norway}
\email{jvrennemo@gmail.com}

\begin{document}
\begin{abstract}
The Grassmannian $\Gr(2,5)$ is embedded in $\PP^{9}$ via the Pl\"ucker embedding.
The intersection of two general $\PGL(10)$-translates of $\Gr(2,5)$ is a Calabi--Yau 3-fold $X$, and the intersection of the projective duals of the two translates is another Calabi--Yau 3-fold $Y$, deformation equivalent to $X$. Applying results of Kuznetsov and Jiang--Leung--Xie shows that $X$ and $Y$ are derived equivalent, which by a result of Addington implies that their third cohomology groups are isomorphic as polarised Hodge structures. We show that $X$ and $Y$ provide counterexamples to a certain ``birational'' Torelli statement for Calabi--Yau 3-folds, namely, they are deformation equivalent, derived equivalent, and have isomorphic Hodge structures, but they are not birational.
\end{abstract}
\maketitle

\section{Introduction}
We study the following pair of Calabi--Yau 3-folds.
Let $V = \CC^{5}$, and consider the Grassmannian $\Gr(2,5)$, contained in $\PP(\wedge^{2}V) = \PP^{9}$ via the Pl\"ucker embedding.
Let $g \in \PGL(\wedge^{2} V)$ be a general element, and let 
\[
X_{g} = \Gr(2,5) \cap g\Gr(2,5).
\]
Then $X_{g}$ is a simply connected smooth Calabi--Yau 3-fold.
The non-trivial Hodge numbers of $X_{g}$ were computed by Kanazawa \cite{kanazawa_pfaffian_2012} to be $h^{1,1}(X) = 1$, $h^{1,2}(X) = 51$.
The family of all $X_{g}$ is locally complete.

Let us write $\Gr(2,5)^{\vee} \subset \PP(\wedge^{2}V^{\vee})$ for the projective dual of $\Gr(2,5)$.
If we choose an isomorphism $V \cong V^{\vee}$ and so identify $\PP(\wedge^{2} V)$ with $\PP(\wedge^{2} V^{\vee})$, we have $\Gr(2,5) = \Gr(2,5)^{\vee}$.
We define 
\[
Y_{g} = \Gr(2,5)^{\vee} \cap (g\Gr(2,5))^{\vee}.
\]
Under the identification of $\PP(\wedge^{2}V)$ with $\PP(\wedge^{2}V^{\vee})$, we have
\[
Y_{g} = \Gr(2,5) \cap g^{-t}\Gr(2,5),
\]
where $g^{-t}$ is the inverse transpose of $g$.
Obviously, $Y_{g}$ is deformation equivalent to $X_{g}$.

The first result of this paper, which we learned from Kuznetsov, is that these varieties are derived equivalent.
The statement is a corollary of Kuznetsov's result \cite{kuznetsov_hyperplane_2006} that $\Gr(2,5)$ and $\Gr(2,5)^{\vee}$ are homological projective duals, together with Jiang--Leung--Xie's general results \cite{jiang_categorical_2017} on intersections and HP duality.\footnote{Similar results have been announced by Kuznetsov and Perry, see the talk \cite{kuznetsov_categorical_nodate} and the forthcoming paper \cite{kuznetsov_categorical_forthcoming}.}
\begin{nprop}
Let $g \in \PGL(\wedge^2 V)$ be such that $X_{g}$ and $Y_{g}$ are of expected dimension.
Then we have an equivalence of derived categories
\[
D^{b}(X_g) \cong D^{b}(Y_g)
\]
\end{nprop}
\begin{proof}
Let $S = \Gr(2,5) \subset \PP(\wedge^{2} V)$, let $S_{g} = g\Gr(2,5) \subset \PP(\wedge^{2}V)$, let $T = \Gr(2,5)^{\vee} \subset \PP(\wedge^{2}V^{\vee})$, and let $T_{g} = (g\Gr(2,5))^{\vee} \subset \PP(\wedge^{2} V^{\vee})$.
Then in the language of \cite{jiang_categorical_2017}, the pairs $(S,T)$ and $(S_{g},T_{g})$ are ``admissible'', by \cite[Lem.~3.4 (1)]{jiang_categorical_2017}.

By \cite[Sec.~6.1]{kuznetsov_hyperplane_2006}, there exist Lefschetz decompositions for $S,T,S_{g},T_{g}$ such that $S$ (resp.~$S_{g}$) is homologically projectively dual to $T$ (resp.~$T_{g}$).
Applying the main theorem of \cite{jiang_categorical_2017} then shows that $D^{b}(X_{g}) = D^{b}(S \cap S_{g}) \cong D^{b}(T \cap T_{g}) = D^{b}(Y_{g})$. 
\end{proof}

A result of Addington then implies that $H^{3}(X_{g},\ZZ) \cong H^{3}(Y_{g},\ZZ)$ as polarised Hodge structures, see Prop.~\ref{thm:DerivedEquivalencePreservesHodge}.

The main technical contribution of this paper is the following:

\begin{notbirational}
For a general choice of $g$, the varieties $X_g$ and $Y_g$ are not birational.
\end{notbirational}
Since $X_{g}$ and $Y_{g}$ have Picard number 1, this reduces to showing that $X_{g}$ and $Y_{g}$ are not projectively equivalent.

\subsection{Counterexamples to Torelli}
Our interest in this example is due to the fact that it is a counterexample to a Torelli-type statement; namely it is to our knowledge the first example of deformation equivalent, non-birational Calabi--Yau 3-folds with equivalent middle Hodge structures.

We now briefly summarise what is known about the general ``Torelli problem'' for Calabi--Yau 3-folds, by which we mean the question of to what extent a member $X$ of a given deformation family of Calabi--Yau 3-folds is determined by the polarised Hodge structure $H^{3}(X,\ZZ)$.
If $H^{3}(X,\ZZ) \cong H^{3}(Y,\ZZ)$ as polarised Hodge structures, we say that $X$ and $Y$ are ``Hodge equivalent''.
There are two natural ways in which one could hope for a Torelli principle to hold: Given deformation equivalent Calabi--Yau 3-folds $X$ and $Y$ which are Hodge equivalent, one can ask if $X$ and $Y$ are isomorphic (``strong Torelli'') or birational (``birational Torelli'').
The terminology is invented by us for the sake of the current discussion.

The first counterexample to strong Torelli was given by \szendroi{} in \cite{szendroi_calabi-yau_2000}, where he produced pairs of Calabi--Yau 3-folds which are deformation equivalent, Hodge equivalent, and birational, but not isomorphic.
These examples are deformations of resolutions of hypersurfaces in certain weighted projective 4-spaces.

A candidate for counterexamples to birational Torelli were constructed by Aspinwall and Morrison in \cite{aspinwall_chiral_1994} and analysed by \szendroi{} in \cite{szendroi_example_2004}. 
The varieties considered are resolutions of finite group quotients of special quintic hypersurfaces in $\PP^4$; they have fundamental group $\ZZ/5\ZZ$.

There is a natural $\ZZ/5\ZZ$-action on the moduli space, and the elements $Y_{1},\ldots, Y_{5}$ of a given orbit have the same universal covering space $Z$.
\szendroi{} shows that $H^{3}(Y_{i},\ZZ)$ is a sub-Hodge structure of index dividing 25 in $H^{3}(Z,\ZZ)$ for all $i$.
In particular, $H^{3}(Y_{i},\ZZ[\frac 1 5]) \cong H^{3}(Y_{j},\ZZ[\frac 1 5])$ for all $i,j$.
They are furthermore pairwise non-isomorphic and conjecturally non-birational.
The Picard rank of $Y_{i}$ is 5, which makes the conjecture hard to verify.

Further potential counterexamples to birational Torelli were given by \cald{} in \cite{caldararu_non-birational_2007}.
Given a Calabi--Yau 3-fold $X$ admitting a genus 1 fibration, the moduli space $X^{(k)}$ of line bundles of degree $k$ on the fibres is a new Calabi--Yau 3-fold, fibred in genus 1 over the same base.
Bridgeland and Maciocia show \cite{bridgeland_fourier-mukai_2002} that the $X^{(k)}$ for different $k$ are all derived equivalent (thus by Prop.~\ref{thm:DerivedEquivalencePreservesHodge} Hodge equivalent), and \cite{caldararu_non-birational_2007} shows that if a few conditions are imposed on $X$, then the $X^{(k)}$ are not all birational.
The question of whether there exist $X$ satisfying these conditions appears to be open.

If we drop the Calabi--Yau condition, then birational Torelli is already known to fail.
Uehara has constructed counterexamples which are minimal 3-folds with Kodaira dimension 1 \cite{uehara_counterexample_2012}.

\subsection{Rational, derived and Hodge}
It's interesting to compare three natural relations between pairs of Calabi--Yau 3-folds $X$ and $Y$:
\begin{enumerate}
\item $X$ and $Y$ are birational
\item $D(X) \cong D(Y)$
\item $H^{3}(X,\ZZ) \cong H^{3}(Y,\ZZ)$ as polarised Hodge structures
\end{enumerate}
By Bridgeland's result \cite{bridgeland_flops_2002}, we know that $(1) \implies (2)$, and by Prop.~\ref{thm:DerivedEquivalencePreservesHodge} below we know that $(2) \implies (3)$.
This argument also gives an alternative proof of the non-trivial implication $(1) \implies (3)$, which was shown by Koll\'ar \cite[Cor.~4.12]{kollar_flops_1989}.

There are several counterexamples to $(2) \implies (1)$: The Pfaffian--Grassmannian pairs \cite{rodland_pfaffian_2000,borisov_pfaffian-grassmannian_2009,kuznetsov_homological_2006}, Hosono and Takagi's example \cite{hosono_double_2016}, and the example studied by Gross--Popescu \cite{gross_calabiyau_2001}, Bak \cite{bak_spectral_2009} and Schnell \cite{schnell2013fundamental}.
The example of this paper is the only one we know of for which the varieties are deformation equivalent.

The implication $(3) \implies (2)$ obviously cannot hold without further restriction, e.g.~since the rank of $H^{\text{even}}$ is a derived invariant.
It is natural to ask if $(3) \implies (2)$ could hold for deformation-equivalent $X$ and $Y$.
The computations of \cite{aspinwall_chiral_1994} seem to suggest that the examples considered there could provide counterexamples.

\subsection{A weaker Torelli statement}
We may compare the state of the Torelli problem for Calabi--Yau 3-folds to that for hyper-K\"ahler varieties.
For a HK variety $X$, the relevant Hodge-theoretic data is $H^{2}(X,\ZZ)$ equipped with the Beauville--Bogomolov--Fujiki quadratic form.
Examples of deformation-equivalent, non-birational pairs of HK varieties with equivalent $H^{2}(X,\ZZ)$ were first discovered by Namikawa \cite{namikawa_counter-example_2002}.
This failure of the Torelli principle can be remedied by Verbitsky's Torelli theorem \cite{verbitsky_mapping_2013,Huybrechts}, which among other things asserts that if $X$ and $Y$ are hyper-K\"ahler and $H^{2}(X,\ZZ) \to H^{2}(Y,\ZZ)$ is a Hodge equivalence induced by parallel transport along some path in the moduli space of varieties, then $X$ and $Y$ are birational.
With a view to a potential similar Torelli theorem for Calabi--Yau 3-folds, it would be interesting to know if our isomorphism $H^{3}(X_{g},\ZZ) \cong H^{3}(Y_{g},\ZZ)$ can be induced by parallel transport along some deformation from $X_{g}$ to $Y_{g}$.

\subsection{Other work}
The 3-folds studied in this paper have received some attention in recent years, partly from the perspective of mirror symmetry, see \cite{kapustka,kapustka_geometric_2011,kanazawa_pfaffian_2012, inoue}.
In particular, the question of whether the $X_{g}$ have non-trivial Fourier--Mukai partners was raised in \cite[Rmk.~p.~25]{hosono_mirror_2014}.

\subsection{Acknowledgements}
We learned of this example from Sasha Kuznetsov, who stated the derived equivalence $D(X_{g}) \cong D(Y_{g})$ and conjectured that $X_{g}$ and $Y_{g}$ are not birational.
Nick Addington pointed out \cite[footnote p.~857]{addington_moduli_2016} to us and observed that this pair would give a new kind of counterexample to CY3 Torelli.
We thank them both for useful conversations, along with Michał Kapustka, Laurent Manivel, Alexander Polishchuk, Kristian Ranestad, Balázs Szendrői and Richard Thomas.

After this paper was completed, we learned that similar results have been obtained by Borisov, \cald{}, and Perry, see \cite{borisov_intersections_2017}.

\section{Derived equivalence implies Hodge equivalence}
\label{sec:derivedImpliesHodge}
\begin{nprop}
\label{thm:DerivedEquivalencePreservesHodge}
Let $X$ and $Y$ be smooth, complex, projective varieties of odd dimension $n$ such that $H^{k}(X,\ZZ) = H^{k}(Y,\ZZ) = 0$ for all odd $k < n$.
If $D^{b}(X) \cong D^{b}(Y)$, then 
\[
H^n(X,\ZZ)/\text{torsion} \cong H^n(Y,\ZZ)/\text{torsion}
\]
as polarised Hodge structures.

In particular, if $X$ and $Y$ are smooth, projective 3-folds with $H^{1}(X,\QQ) = 0$ and $D^{b}(X) \cong D^{b}(Y)$, then 
\[
H^3(X,\ZZ)/\text{torsion} \cong H^3(Y,\ZZ)/\text{torsion}
\]
as polarised Hodge structures.
\end{nprop}
\begin{nremark}
By \cite{popa_derived_2011}, we have $D(X) \cong D(Y) \implies H^{1}(X,\QQ) \cong H^{1}(Y,\QQ)$, which is why we only require $H^{1}(X,\QQ) = 0$ and not $H^{1}(Y,\QQ) = 0$ in the 3-fold case.
\end{nremark}

\begin{nremark}
With the assumptions we take it is not true that the torsion part of $H^{n}$ is preserved; see \cite{addington_brauer_2017}.
However, if we assume that $H^{k}(X,\ZZ) = H^{k}(Y,\ZZ) = 0$ for all odd $k \not= n$, so that in both varieties all torsion lives in $H^{n}$ and $H^{n+1}$, then the proof gives the stronger claim that $H^{n}(X,\ZZ) \cong H^{n}(Y,\ZZ)$.
\end{nremark}
The identification of rational polarised Hodge structures $H^{n}(X,\QQ) \cong H^{n}(Y,\QQ)$ follows from the techniques of \cite[Ch.~5]{huybrechts_fourier-mukai_2006}; the novelty here is that we recover the integral Hodge structure.
The idea of the proof is to observe that under our assumptions the image of $H^{n}(X,\ZZ)$ in $H^{n}(X,\QQ)$ agrees with the image of the Chern character map from (odd-degree) complex $K$-theory; the latter is always preserved by derived equivalence.

We learned the idea of using complex $K$-theory in this way from Nick Addington.
Our proof is along the lines of \cite[footnote p.~857]{addington_moduli_2016}, which treats the case of Calabi--Yau 3-folds.

\subsection{Cohomology and $K$-theory}
Before proving Prop.~\ref{thm:DerivedEquivalencePreservesHodge}, we recall some facts regarding which cohomological structures of a variety are preserved by derived equivalence.

Let $X$ and $Y$ be smooth, projective varieties.
We let $K^{*}(X)$ denote the (even or odd) complex $K$-theory group of the underlying topological space.

For a class $\alpha \in K^{*}(X)$, define the Mukai vector $v(\alpha) \in H^{*}(X,\QQ)$ by
\[
v(\alpha) = \ch(\alpha)\sqrt{\td_{X}} \in H^{*}(X,\QQ).
\]
The square root is defined by treating $\td_{X}$ as a formal power series in indeterminates, requiring that the initial term of $\sqrt{\td_{X}}$ is $1$.

Given an object $\cE \in D^{b}(X \times Y)$, we get a Fourier--Mukai functor
\[
\Phi_{\cE, X \to Y} = (\pi_{Y})_{*}(\cE \otimes \pi_{X}^{*}(-)) \colon D^{b}(X) \to D^{b}(Y).
\]
Representing $\cE$ by locally free sheaves and taking the alternating sum of the complex vector bundles appearing gives a well-defined class $[\cE] \in K^{0}(X)$.
Thus we get a map in $K$-theory\footnote{See e.g.~\cite[Ch.~IV.5]{karoubi_k-theory_1978} for the Gysin map in $K$-theory, which a priori depends on the complex structures of $X$ and $Y$.}
\[
\Phi_{\cE, X \to Y}^{K} = (\pi_{Y})_{!}([\cE] \otimes \pi_{X}^{*}(-)) \colon K^{*}(X) \to K^{*}(Y),
\]
and a map in rational cohomology
\[
\Phi_{\cE, X \to Y}^{H} = (\pi_{Y})_{!}(v(\cE) \cup \pi_{X}^{*}(-)) \colon H^{*}(X,\QQ) \to H^{*}(Y,\QQ).
\]

Using the differentiable Riemann--Roch theorem \cite[Sec.~3]{atiyah_vector_1961}, \cite[V.4.18]{karoubi_k-theory_1978}, the same proof as for \cite[Cor.~5.29]{huybrechts_fourier-mukai_2006} shows that the following square commutes:
\begin{equation}
\label{eqn:MukaiVectorCommutes}
\begin{tikzcd}
K^{*}(X) \arrow[r, "\Phi^{K}_{\cE, X \to Y}"] \arrow[d, "v(-)"] & K^{*}(Y) \arrow[d, "v(-)"]\\
H^{*}(X,\QQ) \arrow[r, "\Phi^{H}_{\cE, X \to Y}"] & H^{*}(Y,\QQ)
\end{tikzcd}
\end{equation}

Define the convolution kernel $\cE \circ \cF$ in the standard way \cite[Prop.~5.10]{huybrechts_fourier-mukai_2006}, so that $\Phi_{\cE \circ \cF, X \to Z} = \Phi_{\cE, Y \to Z} \circ \Phi_{\cF, X \to Y}$.
We have
\begin{equation}
\label{eqn:KFunctoriality}
\Phi^{K}_{\cE \circ \cF, X \to Z} = \Phi^{K}_{\cE, Y \to Z} \circ \Phi^{K}_{\cF, X \to Y},
\end{equation}
\[
\Phi^{H}_{\cE \circ \cF, X \to Z} = \Phi^{H}_{\cE, Y \to Z} \circ \Phi^{H}_{\cF, X \to Y}
\]
and
\[
\Phi^{K}_{\cO_{\Delta}, X \to X} = \id_{K^{*}_{\ttop}(X)}, \ \ \ \Phi^{H}_{\cO_{\Delta}, X \to X} = \id_{H^{*}(X, \QQ)}.
\]
The relation (\ref{eqn:KFunctoriality}) uses the fact that the Gysin map in complex $K$-theory agrees with the push-forward in algebraic $K$-theory, see \cite{atiyah_riemannroch_1962}.

If $\cE \in D(X\times Y)$ determines an equivalence $\Phi_{\cE, X \to Y} \colon D(X) \cong D(Y)$, the inverse is given by $\Phi_{\cF, Y \to X}$, where $\cF = \cE^{\vee} \otimes \pi_{Y}^{*}(\omega_{Y})[\dim Y]$ and $\cF \circ \cE = \cO_{\Delta_{X}}$, $\cE \circ \cF = \cO_{\Delta_{Y}}$.
As a formal consequence, we have 
\begin{nlemma}
 If $\Phi_{\cE, X \to Y}$ is an equivalence, then $\Phi^{K}_{\cE, X \to Y}$ and $\Phi^{H}_{\cE, X \to Y}$ are isomorphisms.
\end{nlemma}

\subsection{The generalised Mukai pairing}
Following \cald{} \cite{caldararu_mukai_2005}, we may define a generalised Mukai pairing on $H^{*}(X,\CC)$, in the following way.
If $\alpha \in H^{d}(X,\CC)$, let $\alpha^{\vee} = \sqrt{-1}^{d}\alpha$, and extend this operation linearly to get an automorphism $(-)^{\vee}$ of $H^{*}(X,\CC)$.
Given classes $\alpha, \beta \in H^{*}(X,\CC)$, define the Mukai pairing by
\[
(\alpha, \beta) = \int_{X} e^{c_{1}(X)/2} \cup \alpha^{\vee} \cup \beta.
\]

\begin{nlemma}[{\cite[Prop.~5.39, Prop.~5.44]{huybrechts_fourier-mukai_2006}}]
If $\Phi_{\cE,X \to Y}$ is an equivalence, then the map $\Phi^{H}_{\cE,X\to Y} \colon H^{*}(X,\CC) \to H^{*}(Y,\CC)$ is an isometry with respect to the Mukai pairing, and for each $k$ takes the subspace $\bigoplus_{p-q = k} H^{p,q}(X)$ isomorphically to $\bigoplus_{p-q = k}H^{p,q}(Y)$.
\end{nlemma}

\subsection{The image of the Chern character}
Let $X$ be a finite CW complex, and let $X_{\le d}$ denote the $d$-skeleton of $X$.
There is a filtration on $K^{*}(X)$ defined by letting the subspace $K^{*}_{>d}(X) = K^{*}_{\ge d+1}(X)$ be the kernel of the restriction map $K^{*}(X) \to K^{*}(X_{\le d})$.
The associated graded group is computed by the Atiyah--Hirzebruch spectral sequence \cite{atiyah_vector_1961}, which has $E_{2}$-page given by $E^{p,q}_{2} = H^{p}(X,K^{q}(\pt))$, and $E^{p,q}_{\infty} = K_{\ge p}^{p+q}(X)/K_{\ge p+1}^{p+q}(X)$.

Let $d_{r}$ be the differential on the $E_{r}$-page.
By \cite{atiyah_vector_1961}, we have $d_{r} \otimes \QQ = 0$ for $r\ge 2$, and so $E^{p,q}_{\infty}$ is a subquotient of $E^{p,q}_{2}$ which becomes isomorphic after tensoring with $\QQ$.
We thus get a well-defined map $I \colon E^{p,0}_{\infty} \to E^{p,0}_{\infty}\otimes \QQ = E^{p,0}_{2} \otimes \QQ = H^{p}(X,\QQ)$.
Let $\ch_{p}$ denote the piece of the Chern character which lands in $H^{p}(X,\QQ)$, and let $\iota \colon H^{*}(X,\ZZ) \to H^{*}(X,\QQ)$ be the scalar extension map.
\begin{nlemma}
\label{thm:ChIsSurjective}
The image of $\ch_{p} \colon K_{\ge p}^{p}(X) \to H^{p}(X,\QQ)$ equals $I(E^{p,0}_{\infty})$.
\end{nlemma}
\begin{proof}
The proof is essentially the same as that of \cite[Cor.~2.5]{atiyah_vector_1961}, which says that if $H^{*}(X,\ZZ)$ is torsion-free, then the image of $\ch_{p}$ is $H^{p}(X,\ZZ)$.

Filtering $X$ by $p$-skeletons gives a trivial spectral sequence $F^{p,q}_{i}$ converging to $H^{p+q}(X,\QQ)$ with $F^{p,0}_{2} = H^{p}(X,\QQ)$, $F^{p,q}_{2} = 0$ for $q \not= 0$, and degenerating at the $F_{2}$-page.
The Chern character defines a map of spectral sequences $E^{p,*}_{i} \to F^{p,*}_{i}$, whose image in $F^{p,0}_{2} = H^{p}(X,\QQ)$ is the same as $\iota(H^{p}(X,\ZZ))$.
We then get
\[
\ch_{p}(K_{\ge p}^{p}(X)) = \ch(E^{p,0}_{\infty}) = I(E^{p,0}_{\infty}) \subset H^{p}(X,\QQ) = F^{p,0}_{\infty}.
\]
\end{proof}

\begin{nlemma}
\label{thm:imageOfChIsIntegralStructure}
Let $X$ be a smooth, projective variety of odd dimension $n$ such that $H^{k}(X,\ZZ) = 0$ for all odd $k < n$. Then
\[
v(K^{1}(X)) = \ch(K^{1}(X)) = \iota(H^{n}(X,\ZZ)).
\]
\end{nlemma}
\begin{proof}
The first equality is obvious from the definition, since $H^{\odd}(X,\QQ) = H^{n}(X,\QQ)$.

Note that $K^{1}(X) = K^{1}_{\ge n}(X)$, since $H^{\odd}(X_{\le n-1},\ZZ) = 0$, and so $K^{1}(X_{\le n-1}) = 0$.
Thus, by Lemma \ref{thm:ChIsSurjective}, we have $\ch_{n}(K^{1}(X)) = \ch_{n}(K^{n}(X)) = I(E^{n,0}_{\infty}(X))$.

By Poincar\'e duality and the universal coefficient theorem, we have that $H^{k}(X,\ZZ)$ is torsion free for any even $k \ge n + 3$.
This implies that $d_{k}|_{E^{n,0}_{k}} = 0$ for all $k \ge 2$.
Hence $E^{n,0}_{\infty} = H^{n}(X,\ZZ)$ and so $I(E^{n,0}_{\infty}) = \iota(H^{n}(X,\ZZ))$.
\end{proof}

\begin{proof}[Proof of Prop.~\ref{thm:DerivedEquivalencePreservesHodge}]
By Orlov's representability theorem \cite[Thm.~3.2.1]{orlov_derived_2003}, the equivalence $D^{b}(X) \to D^{b}(Y)$ is induced by a Fourier--Mukai kernel $\cE \in D^{b}(X\times Y)$.

Thus $\Phi^{H}_{\cE, X \to Y}$ gives an isomorphism $H^{*}(X,\QQ) \to H^{*}(Y,\QQ)$ which preserves the generalised Mukai pairing and the subspaces $H_{k}^{*} = \oplus_{p - q = k} H^{p,q}$.
As $H^{\odd}(-,\QQ) = H^{n}(-,\QQ)$ for both $X$ and $Y$, this implies that $\Phi_{\cE, X \to Y}^{H}$ identifies $H^{n}(X,\QQ)$ with $H^{n}(Y,\QQ)$.
The Mukai pairing on $H^{n}(X,\CC)$ and $H^{n}(Y,\CC)$ is a rescaling of the standard inner product, and the decomposition $H^{n} = \oplus_{k} H_{k}^{*}$ is the Hodge decomposition, so that $\Phi^{H}_{\cE, X\to Y}$ identifies $H^{n}(X,\QQ)$ and $H^{n}(Y,\QQ)$ as polarised rational Hodge structures.

It remains to see that $\Phi^{H}_{\cE,X\to Y}(\iota(H^{n}(X,\ZZ))) = \iota(H^{n}(Y,\ZZ))$.
This follows from Lemma \ref{thm:imageOfChIsIntegralStructure}, the equality $\Phi^{K}_{\cE, X \to Y}(K^{1}(X)) = K^{1}(Y)$, and the commutativity of (\ref{eqn:MukaiVectorCommutes}).
\end{proof}

\section{Some observations on $X_{g}$}
Let $\PGL(\wedge^{2}V)^{\circ} \subset \PGL(\wedge^{2}V)$ be the open subset of $g$ such that $X_{g}$ is non-singular, and let $\cX \to \PGL(\wedge^{2}V)^{\circ}$ be the family of all such $X_{g}$.
\begin{nprop}
  The family $\cX \to \PGL(\wedge^{2}V)^{\circ}$ is locally complete.
\end{nprop}
\begin{proof}
For any point $g \in \PGL(\wedge^{2}V)^{\circ}$, we must show that the Kodaira--Spencer map $T_{\PGL, g} \to H^{1}(T_{X_{g}})$ is surjective.
We may identify $T_{\PGL, g}$ with $H^{0}(\PP^9,T_{\PP^9})$ in such a way that the Kodaira--Spencer map becomes the composition 
\[
H^{0}({\PP^9},T_{{\PP^9}}) \to H^{0}(\Gr(2,5),N_{g\Gr|{\PP^9}}) \to H^{0}(X_{g},N_{X_{g}|\Gr}) \to H^{1}(X_{g},T_{X_{g}}),
\]
where the first two maps are the natural restriction maps, and the last is the boundary map in a long exact sequence.
The map $H^{0}(T_{{\PP^9}}) \to H^{0}(N_{X_{g}|\Gr})$ factors as
\[
H^{0}(\PP^{9},T_{{\PP^{9}}}) \stackrel{\phi_{1}}{\to} H^{0}(X_{g},T_{{\PP^{9}}}|_{X_{g}}) \stackrel{\phi_{2}}{\to} H^{0}(X_{g},N_{g\Gr|{\PP^{9}}}|_{X_{g}}) \stackrel{\phi_{3}}{\to} H^{0}(X_{g},N_{X_{g}|\Gr}).
\]
Here $\phi_{1}$ and $\phi_{2}$ are surjective by Lemma \ref{thm:someVanishing} below, and the sheaves $N_{g\Gr|{\PP^{9}}}|_{X_{g}}$ and $N_{X_{g}|\Gr}$ are isomorphic, so $\phi_{3}$ is an isomorphism.
Finally the map $H^{0}(N_{X_{g}|\Gr}) \to H^{1}(T_{X_{g}})$ is surjective by Lemma \ref{thm:someVanishing}.
\end{proof}

\begin{nlemma}
\label{thm:someVanishing}
We have $H^{1}(X, T_{\Gr(2,5)}|_{X}) = 0$.
The map $H^{0}(T_{{\PP^{9}}}) \to H^{0}(T_{{\PP^{9}}}|_{X})$ is surjective.
\end{nlemma}
\begin{proof}
The sheaf $\cO_{\Gr(2,5)}$ on $\PP^{9}$ admits a Pfaffian resolution
\begin{equation}
\label{eqn:PfaffianResolution}
0 \to \cO_{\PP^{9}}(-5) \to \cO_{\PP^{9}}(-3)^{\oplus 5} \to \cO_{\PP^{9}}(-2)^{\oplus 5} \to \cO_{\PP^{9}} \to \cO_{\Gr(2,5)} \to 0,
\end{equation}
which gives a resolution of $T_{\Gr(2,5)}|_{X}$ on $\Gr(2,5)$ as follows:
\[
0 \to T_{\Gr}(-5) \to T_{\Gr}(-3)^{\oplus 5} \to T_{\Gr}(-2)^{\oplus 5} \to T_{\Gr} \to T_{\Gr}|_{X} \to 0.
\]
Using the Borel--Weil--Bott theorem for vector bundles on $\Gr(2,5)$ (as in \cite{kuznetsov_homological_2006}), we find that 
\[
H^i(T_{\Gr(2,5)}{(-t)})=0
\] for all $t\ge 0$ and $1\le i\le 4$ and $H^5(T_{\Gr(2,5)}(-5))=\CC$. The hypercohomology spectral sequence for the resolution of $T_{\Gr}|_{X}$ above then shows that $H^{1}(X,T_{\Gr}|_{X}) = 0$.

Since the intersection of the two Grassmannians is transversal, the structure sheaf $\cO_{X}$ admits a resolution on $\PP^{9}$ obtained by taking the tensor product of resolutions (\ref{eqn:PfaffianResolution}) for $\cO_{\Gr(2,5)}$ and $\cO_{g\Gr(2,5)}$.
This gives a resolution of $T_{{\PP^{9}}}|_{X}$ on $\PP^{9}$ which looks like
\[
0 \to T_{\PP^{9}}(-10) \to \cdots \to T_{\PP^{9}}(-2)^{\oplus 10} \to T_{\PP^{9}} \to T_{\PP^{9}}|_{X} \to 0.
\]
Using this resolution, a similar argument to the one above gives the second claim.
\end{proof}

\begin{nlemma}
\label{thm:topologicalInvariants}
The variety $X_{g}$ is simply connected, $\Pic(\Gr(2,5))\to \Pic(X_{g})$ is an isomorphism, and $H^{3}(X_{g},\ZZ)$ is torsion free.
\end{nlemma}
\begin{proof}
We use a result of Sommese \cite[Cor.~on first page]{sommese_complex_1982}, which implies that for subvarieties $A,B$ of projective space, the relative homotopy groups $\pi_j(A,A\cap B)$ vanish for $j\le \min\{\dim A,\dim B+1\}-\codim B$. In our case, when $A$ and $B$ are general $\Gr(2,5)$-translates in $\PP^9$, this gives
\[
\pi_{\le 3}(\Gr(2,5),X_{g}) = 0.
\] 
In particular, $\pi_{1}(X_{g}) \to \pi_{1}(\Gr(2,5))$ is an isomorphism.
By the relative Hurewicz theorem, we have that $H_{\le 3}(\Gr(2,5),X_{g}) = 0$.
The remaining claims now follow from the universal coefficient theorem and the cohomology long exact sequence of the pair $(\Gr(2,5),X_{g})$.
\end{proof}

\section{$X_{g}$ and $Y_{g}$ are not birational}
\label{notbirational}
Let $X_{g}$ and $Y_{g}$ be defined as in the introduction.
We choose a basis $\{e_{i}\}_{i=1}^{5}$ for $V$, which gives a dual basis $\{e^{i}\}_{i=1}^{5}$ for $V^{\vee}$.
This choice gives us a natural isomorphism 
\begin{equation}
\label{eqn:identification}
\PP(\wedge^{2} V^{\vee}) \cong \PP(\wedge^{2}V).
\end{equation}
We let $g^{-t} \in \PGL(\wedge^{2}V)$ be the inverse transpose of $g$ with respect to the basis $\{e_{i}\wedge e_{j}\}$.
Under the isomorphism (\ref{eqn:identification}), the action of $g$ on $\PGL(\wedge^{2}V^{\vee})$ is identified with the action of $g^{-t}$ on $\PGL(\wedge^{2}V)$.
Thus under (\ref{eqn:identification}), we have
\[
Y_{g} = \Gr(2,5) \cap g^{-t}\Gr(2,5) \subset \PP(\wedge^{2}V).
\]
From this point on we will always think of $Y_{g}$ as embedded in $\PP(\wedge^{2}V)$ in this way.

The purpose of this section is to show
\begin{nthm}
\label{thm:XAndYNotBirational}
For a general choice of $g$, the varieties $X_{g}$ and $Y_{g}$ are not birational.
\end{nthm}
\begin{proof}
By Lemma \ref{thm:topologicalInvariants}, $\Pic(X_g)\simeq \Pic(Y_g)\simeq \ZZ H$ where $H$ is the hyperplane section in the Pl\"ucker embedding. If $X_g$ and $Y_g$ are birational, they are isomorphic in codimension 1 (since they are minimal models). 
Now $H$ is very ample on both $X_g$ and $Y_g$, so $X_g$ and $Y_g$ are isomorphic.
Since this isomorphism induces an isomorphism on $H^0(X_g,\cO(H))$, it is induced from a linear automorphism of $\PP(\wedge^{2}V)$. 

Hence we reduce to showing that $X_{g}$ and $Y_{g}$ are not projectively equivalent. Proposition \ref{thm:XIsContainedInOnlyTwoGrassmannians} shows that $X_{g}$ and $Y_{g}$ are each contained in exactly 2 translates of $\Gr(2,5)$, hence if $X_{g}$ and $Y_{g}$ are projectively equivalent, then there must be a projective transformation of $\PP(\wedge^{2} V)$ identifying the unordered pairs $(\Gr(2,5), g\Gr(2,5))$ and $(\Gr(2,5), g^{-t}\Gr(2,5))$. 
Finally, Lemma \ref{thm:GenericallyProjectivelyDifferent} shows that no such transformation exists.
\end{proof}

From this point on, we fix a general $g \in \PGL(\wedge^{2}V)$ and write $X = X_{g}$, $Y = Y_{g}$.

\subsection{$X$ is contained in only two Grassmannian translates}
We consider $|\cO(2)| = \PP(H^0(\PP(\wedge^2 V), \cO(2))$, the space of quadrics in $\PP (\wedge^2 V)$, and 
\[
|\cI_{X}(2)|, |\cI_{\Gr}(2)|, |\cI_{g\Gr}(2)| \subset |\cO(2)|,
\]
the linear systems of quadrics containing $X$, $\Gr(2,5)$ and $g\Gr(2,5)$ respectively.

\begin{nprop}
\label{thm:XIsContainedInOnlyTwoGrassmannians}
If there exists an $h \in \PGL(\wedge^2 V)$ such that $X \subset h\Gr(2,5)$, then either $h\Gr(2,5) = \Gr(2,5)$ or $h\Gr(2,5) = g\Gr(2,5)$.
\end{nprop}
Let's briefly summarise the proof. We first observe that $|\cI_{\Gr}(2)| \cong |\cI_{g\Gr}(2)| \cong \PP^{4}$, that $\Span(|\cI_{\Gr}(2)|, |\cI_{g\Gr}(2)|) = |\cI_{X}(2)| \cong \PP^{9}$, and that all quadrics in $|\cI_{\Gr}(2)|$ and $|\cI_{g\Gr}(2)|$ have rank 6.

We assume for a contradiction that $X$ is contained in a third translate $h\Gr(2,5)$, which means that $|\cI_{h\Gr}(2)| \subset |\cI_{X}(2)|$.
A general point of $|\cI_{h\Gr}(2)|$ lies on a unique line between $|\cI_{\Gr}(2)|$ and $|\cI_{g\Gr}(2)|$, and this defines a rational map $\varphi \colon |\cI_{h\Gr}(2)| \dashrightarrow |\cI_{\Gr}(2)| \times |\cI_{g\Gr}(2)|$.
Lemma \ref{thm:CommonSingularPoint} shows that a pencil of quadrics containing 3 quadrics of sufficiently low rank must have a point at which all the quadrics are singular; in particular any line which intersects all three $|\cI_{\Gr}(2)|, |\cI_{g\Gr}(2)|$ and $|\cI_{h\Gr}(2)|$ is of this kind.

Hence $\overline{\varphi(\cI_{h\Gr}(2))} \subset S \subset |\cI_{\Gr}(2)| \times |\cI_{g\Gr}(2)|$, where $S$ is the correspondence of quadrics with a common singular point.
We can describe the geometry of $S$ rather explicitly, and it's not too hard to derive a contradiction from this point.

\subsubsection{The quadrics defining $\Gr(2,5)$}
We begin with some classical observations about the quadrics containing $\Gr(2,5)$.
Recall the Pl\"ucker relations for $\Gr(2,5)$:
\begin{nlemma}
\label{thm:quadricsDefiningGr}
There is an isomorphism $\PP V \cong |\cI_{\Gr}(2)|$, defined by sending $[v] \in \PP V$ to the quadric $q_{v}$,
\[
q_{v}(\alpha) = \alpha \wedge \alpha \wedge v, \ \ \ \ \alpha \in \wedge^{2} V.
\]
The quadric $q_{v}$ is of rank 6 in 10 variables, and is singular along $\PP(v \wedge V) \subset \Gr(2,5)$.
\end{nlemma}

So to each quadric $q_{v} \in |\cI_{\Gr}(2)|$, we can assign its set of singular points $\PP^3 \cong \PP(v \wedge V) \subset \Gr(2,5)$.
Let $C \subset |\cI_{\Gr}(2)| \times \Gr(2,5)$ be the associated correspondence.
Let $\pi_{1}$ and $\pi_{2}$ denote the projections from $|\cI_{\Gr}(2)| \times \Gr(2,5)$ to the first and second factor respectively.

Given a closed subvariety $Z \subset |\cI_{\Gr}(2)|$, let 
\[
\psi(Z) = \pi_{2}(C \cap \pi_{1}^{-1}(Z)) \subset \Gr(2,5)
\]
\begin{nlemma}
\label{thm:imageOfCorrespondence}
Let $L \subset |\cI_{\Gr}(2)|$ be a linear subspace.
If $\dim L = 2$, then $\psi(L) = \Gr(2,5) \cap H$ for some hyperplane $H \subset \PP(\wedge^2 V)$.
If $\dim L \ge 3$, then $\psi(L) = \Gr(2,5)$. 
\end{nlemma}
\begin{proof}
This is a direct computation.
Identify $V \cong H^{0}(\PP(\wedge^{2}V), \cI_{\Gr}(2))$ as above.
If $\dim L = 2$, let $v_{1}, v_{2}, v_{3} \in V$ be a basis for the vector space corresponding to $L$, and extend this by $v_{4}, v_{5}$ to get a basis for $V$.
Let $\alpha = \sum_{1 \le i < j \le 5} \alpha^{ij} v_{i} \wedge v_{j} \in \wedge^{2}V$ be of rank 2.
Then we may write $\alpha = \alpha_{1} \wedge \alpha_{2}$, with $\alpha_{i} = \sum_{j = 1}^{5} \alpha_{i}^{j}v_{j}$.

Note that $[\alpha] \in \psi(L)$ if and only if $\alpha \in \langle v_{1}, v_{2}, v_{3} \rangle \wedge V$.
This happens if and only if the matrix $\begin{pmatrix} \alpha_{1}^{4} & \alpha_{1}^{5} \\ \alpha^{4}_{2} & \alpha^{5}_{2}\end{pmatrix}$ is singular, which is if and only if $\alpha^{45} = 0$.

The case $\dim L \ge 3$ follows.
\end{proof}

\begin{nlemma}
We have $|\cI_{\Gr}(2)| \cap |\cI_{g\Gr}(2)| = \varnothing$, and
\[
|\cI_{X}(2)| = \Span( |\cI_{\Gr}(2)|, |\cI_{g\Gr}(2)| ) \cong \PP^{9} \subset |\cO_{\PP(\wedge^{2}V)}(2)|.
\]
\end{nlemma}
\begin{proof}
We first show that $|\cI_{\Gr}(2)| \cap |\cI_{g\Gr}(2)| = \varnothing$:
If $q \in |\cI_{\Gr}(2)| \cap |\cI_{g\Gr}(2)|$, then $\psi(q)$ is contained in $\Gr(2,5)$ and $g\Gr(2,5)$, and so we have $\psi(q) = \PP^{3} \subseteq X$, which is impossible.

It remains to see that $|\cI_{X}(2)| = \Span(|\cI_{\Gr}(2)|, |\cI_{g\Gr}(2)|)$.
The inclusion $\supseteq$ is obvious, and to get the inclusion $\subseteq$, it is enough to show $\dim |\cI_{X}(2)| = 9$. Now, the ideal sheaf of each Grassmannian has a Pfaffian resolution of the form
$$
0\to \mathcal O(-5)\to \mathcal O(-3)^5\to \mathcal O(-2)^5\to \cI_{\Gr}\to 0
$$
Since the intersection $X=\Gr(2,5)\cap g \Gr(2,5)$ is transversal, we may restrict the sequence for $g\Gr(2,5)$ above to $\Gr(2,5)$ and so obtain a resolution for $\cI_{X|\Gr}$. 
The Kodaira vanishing theorem shows that $H^{*}(\Gr(2,5),\cO(i)) = 0$ for $-4\le i \le -1$, and so we find that $H^0(\cI_{X|\Gr}(2))=\CC^5$. 
In particular, $H^0(\cI_{X|\PP^9}(2))$ is 10-dimensional, by the sequence
$$
0\to \cI_{\Gr|\PP^9}(2)\to \cI_{X|\PP^9}(2)\to \cI_{X|\Gr}(2)\to 0
$$\end{proof}

\subsubsection{Special pencils of quadrics}
\begin{nlemma}
\label{thm:CommonSingularPoint}
If a pencil of quadrics in $n$ variables contains 3 quadrics of $q_{1},q_{2},q_{3}$ of coranks $c_{1}, c_{2}, c_{3}$ with $c_{1} + c_{2} + c_{3} > n$, then all quadrics in the pencil have a singular point in $\PP^{n-1}$ in common.
\end{nlemma}
\begin{proof}
Let $W = \CC^{n}$, so the quadrics define subsets of $\PP W$.
The pencil of quadrics gives a map of vector bundles on $\PP^{1}$, $\gamma \colon W \otimes \cO_{\PP^{1}}(-1) \to W^{\vee} \otimes \cO_{\PP^{1}}$.
We extend this to an exact sequence
\[
0 \to \cF \to W \otimes \cO_{\PP^{1}}(-1) \to W^{\vee} \otimes \cO_{\PP^{1}} \to \cE \oplus \cT \to 0,
\]
where $\cE$ and $\cT$ are the torsion free and torsion parts of $\mbox{coker }\gamma$, respectively.
Since the map $\gamma$ is defined by a pencil of symmetric forms, we have $\gamma^{\vee}(1) \cong \gamma$, and it follows that $\cF = \cE^{\vee}(-1)$.

Let now $r(-)$ and $d(-)$ denote the rank and degree of a sheaf.
These invariants factor through the $K$-group, so using the above exact sequence gives
\begin{align}
\label{eqn:KTheoryEquality}
\dim W &= d(W^{\vee}\otimes \cO_{\PP^{1}}) - d(W\otimes \cO_{\PP^{1}}(-1)) \nonumber \\
&= d(\cE) - d(\cE^{\vee}(-1)) + d(\cT) = 2d(\cE) + r(\cE) + d(\cT).
\end{align}

We note that $r(\cE)$ is the generic corank of a quadric in the pencil, and so we have $r(\cE) \le \min(c_{1},c_{2},c_{3})$.
At the point $q_{i}$, we have $\dim (\cE|_{q_{i}} \oplus \cT|_{q_{i}}) = c_{i}$, and so $\dim \cT|_{q_{i}} = c_{i} - r(\cE)$.
It follows that $d(\cT) \ge c_{1} + c_{2} + c_{3} - 3r(\cE) > \dim W - 3r(\cE)$.
Combining this with (\ref{eqn:KTheoryEquality}) gives
\[
d(\cE) <  r(\cE).
\]

Now $\cE$ splits into a sum of line bundles $\oplus_{i=1}^{r(\cE)}\cO(d_{i})$, and since $\cE$ is a quotient of a trivial bundle, we have $d_{i} \ge 0$ for all $i$.
The inequality $d(\cE) < r(\cE)$ shows that there must be at least one $i$ such that $d_{i} = 0$.
Correspondingly, the line bundle decomposition of $\cF = \cE^{\vee}(-1)$ has at least one factor $\cO_{\PP^{1}}(-1)$.
The image of this factor inside $W \otimes \cO_{\PP^{1}}(-1)$ must be of the form $w \otimes \cO_{\PP^{1}}(-1)$ for some $w \in W$, which means that $[w] \in \PP W$ is a common singular point for all quadrics in the pencil.
\end{proof}

\begin{proof}[Proof of Prop.~\ref{thm:XIsContainedInOnlyTwoGrassmannians}]
Assume for a contradiction that there exists an $h$ such that $X \subset h\Gr(2,5)$ and $\Gr(2,5) \not= h\Gr(2,5) \not= g\Gr(2,5)$.
Let us write $L = |\cI_{\Gr}(2)|$, $L_{g} = |\cI_{g\Gr}(2)|$ and $L_{h} = |\cI_{h\Gr}(2)|$.
These are all 4-dimensional subspaces of $\PP^{9} \cong |\cI_{X}(2)| = \Span (L, L_{g})$ and we have $L \not= L_{h} \not= L_{g}$.

For any $q \in L_{h} \setminus (L \cup L_{g})$, there exist a unique point $(q_{1},q_{2}) \in L \times L_{g}$ such that $q_{1},q_{2},q$ lie on a line.
This gives a rational map $\varphi \colon L_{h} \dashrightarrow L \times L_{g}$, which is resolved by blowing up the disjoint linear loci $L_{h} \cap L$ and $L_{h} \cap L_{g}$. 

By Lemma \ref{thm:CommonSingularPoint}, the quadrics on the line $\langle q_{1}, q_{2} \rangle$ must all have a singular point in common.
If $x$ is a singular point of $q_{1}$ and $q_{2}$, then as $\Sing(q_{1}) \subset \Gr(2,5)$ and $\Sing(q_{2}) \subset g\Gr(2,5)$, we have $x \in \Gr(2,5) \cap g\Gr(2,5) = X$.

Define the correspondence $\widetilde S \subset L \times L_{g} \times X$ by
\[
\widetilde S = \{(q_{1}, q_{2}, x) \mid x \text{ is a singular point of }q_{1} \text{ and }q_{2}\}
\]
For any $x \in X$, the set of quadrics in $L$ which are singular in $x$ is a line, and similarly the set of quadrics in $L_{g}$ which are singular in $x$ is a line. It follows that $\pi_{X} \colon \widetilde S \to X$ is a $\PP^{1} \times \PP^{1}$-bundle.

Let $S$ be the image of $\widetilde S$ in $L \times L_{g}$.
The above discussion shows that $\varphi(L_{h}) \subset S$.
We thus have the following diagram:
\[
\begin{tikzcd}
& & X \\
L_{h} \arrow[r, dashed, "\widetilde{\varphi}"] \arrow[dr, dashed, "\varphi"]& \widetilde S \arrow[r, hook] \arrow[d, "\pi_{S}"] \arrow[ru, "\pi_{X}"] & L \times L_{g} \times X \arrow[u] \arrow[d]\\
 & S \arrow[r, hook] & L \times L_{g}
\end{tikzcd}
\]
For any rational map $f$, we let $f(L_{h})$ denote the closure of the image of $L_{h}$.

\emph{Claim:} The morphism $\pi_{S}$ is birational, and contracts a finite number of curves.

To see this, note that if $(q_{1},q_{2},x) \not= (q_{1},q_{2},y)$ are contained in $\widetilde S$, then $x, y \in \Sing(q_{1}) \cap \Sing(q_{2})$; in particular $x, y$ lie on a line.
The dimension of the space of lines in $\Gr(2,5)$ is 8, and it follows that since $g$ is chosen general, there are only finitely many lines on $X$.
Any line on $\Gr(2,5)$ is contained in the singular locus of a unique quadric, and so for every line $l$ on $X$ there is a unique pair $(q_{1}, q_{2}) \in L \times L_{g}$ such that $\{q_{1}\} \times \{q_{2}\} \times l$ is contracted by $\pi_{S}$.
This proves the claim.

By construction, the general fibre of $\varphi$ is either a point or a line, and so $\dim \varphi(L_{h}) = \dim \widetilde\varphi(L_{h}) \in \{3,4\}$.
Since the fibres of $\pi_{X}$ are 2-dimensional, we have $\dim \pi_{X}(\widetilde{\varphi}(L_{h})) \in [\dim \widetilde \varphi(L_{h}) - 2, 3]$.
We will show that each of these possible values for the dimensions leads to a contradiction.

\emph{Case 1}, $\dim \varphi(L_{h}) =3$: This means that the general fibre of $\varphi$ is a line.
It follows that $L_{h}$ is the span of two linear subspaces $M \subset L$ and $M_{g} \subset L_{g}$, with $\dim M + \dim M_{g} = 3$.
We may assume that $\dim M > \dim M_{g}$.

If $\dim M = 3$, then using Lemma \ref{thm:imageOfCorrespondence}, we have $\psi(M) = \Gr(2,5) \subseteq h\Gr(2,5)$, and so $\Gr(2,5) = h\Gr(2,5)$, which is a contradiction.

If $\dim M = 2$, then using Lemma \ref{thm:imageOfCorrespondence}, we have $\psi(M) = H \cap \Gr(2,5) \subseteq \Gr(2,5) \cap h\Gr(2,5)$ for some hyperplane $H$.
If $Q$ is a quadric containing $h\Gr(2,5)$, it must therefore contain $H \cap \Gr(2,5)$, and hence $Q \cap \Gr(2,5) = (H \cup H_{Q}) \cap \Gr(2,5)$ for some hyperplane $H_{Q}$.
Choosing 5 general such quadrics $Q_{1}, \ldots, Q_{5}$, we get $X \subset h\Gr(2,5) \cap \Gr(2,5) \subset (H \cup H_{Q_{1}} \cup \cdots \cup H_{Q_{5}}) \cap \Gr(2,5)$.
But $X$ is not contained in any hyperplane, so this is a contradiction.

\emph{Case 2}, $\dim \varphi(L_{h}) = 4, \dim \pi_{X}(\widetilde{\varphi}(L_{h})) = 3$: Since $L_{h}$ is rational and $X$ is Calabi--Yau, this is impossible.

\emph{Case 3}, $\dim \varphi(L_{h}) = 4, \dim \pi_{X}(\widetilde{\varphi}(L_{h})) = 2$: The fibres of $\pi_{X}$ are 2-dimensional, hence 
\begin{equation}
\label{eqn:PushPullEqn}
\widetilde{\varphi}(L_{h}) = \pi_{X}^{-1}\pi_{X}(\widetilde{\varphi}(L_{h})).
\end{equation}
By adjunction, any smooth divisor in $X$ is of general type.
Hence, since $\pi_{X}(\widetilde{\varphi}(L_{h}))$ is unirational, it must be singular, and it then follows from (\ref{eqn:PushPullEqn}) that $\dim \Sing(\widetilde{\varphi}(L_{h})) \ge 2$.
Since $\varphi(L_{h}) = \pi_{S}(\widetilde \varphi(L_{h}))$ and $\pi_{S}$ is an isomorphism outside of finitely many curves, this means that $\varphi(L_{h})$ is also singular.
But $\varphi(L_{h})$ is a blow-up of $L_{h}$ in two disjoint linear subspaces, hence nonsingular.
\end{proof}

\subsection{No linear transformation preserves the pair of Grassmannians}
Recall that we have chosen a basis $\{e_{i}\}_{i=1}^{5}$ of $V$, which gives a dual basis $\{e^{i}\}$ of $V^{\vee}$.
We let $e_{ij} = e_{i} \wedge e_{j}$ and $e^{ij} = e^{i} \wedge e^{j}$ denote the dual (up to scale) basis vectors for $\wedge^{2}V$ and $\wedge^{2} V^{\vee}$.

\def\GG{\mathrm{Gr}}
\begin{nlemma}
\label{thm:GenericallyProjectivelyDifferent}
Let $g \in \PGL(\wedge^{2}V)$ be general.
Then there exists no $h \in\PGL(\wedge^{2}V)$ such that 
\begin{equation}
\label{eqn:GrIsPreserved}
(h\GG(2,5), hg\GG(2,5)) = (\GG(2,5), g^{-t}\GG(2,5))
\end{equation}
or 
\begin{equation}
\label{eqn:GrIsNotPreserved}
(h\GG(2,5), hg\GG(2,5)) = (g^{-t}\GG(2,5), \GG(2,5)).
\end{equation}
\end{nlemma}
\begin{proof}
Consider the action of $\PGL(V) \times\PGL(V)^{\op}$ on $\PGL(\wedge^{2}V)$ given by
\[
((h,h'),g) \mapsto hgh'\ \ \ \ \ h,h' \in\PGL(V), g \in\PGL(\wedge^{2}V).
\]

An element of $\PGL(\wedge^{2}V)$ preserves $\Gr(2,5)$ if and only if it is contained in $\PGL(V) \subset\PGL(\wedge^{2}(V))$.
It follows that the existence of a $h$ satisfying (\ref{eqn:GrIsPreserved}) (resp.~(\ref{eqn:GrIsNotPreserved})) is equivalent to $g$ being in the same $\PGL(V) \times\PGL(V)^{\op}$-orbit as $g^{-t}$ (resp.~$g^{t}$).
So our goal is to show that $g$ is in a different $\PGL(V) \times \PGL(V)^{\op}$-orbit to both $g^{t}$ and $g^{-t}$.

We'll first prove the analogous claim with $\SL$ instead of $\PGL$; the argument for passing to $\PGL$ is given at the end of the proof.
So let $G = \SL(V) \times \SL(V)^{\op}$, which acts on $\SL(\wedge^{2}V)$ as above, and let $g \in \SL(\wedge^{2}V)$ be general.
We aim to show that $g$ is in a different $G$-orbit to both $g^{t}$ and $g^{-t}$.

\emph{Proof that $g$ is not in the $G$-orbit of $g^{-t}$:}
It is enough to construct a function $f \in \CC[\SL(\wedge^{2}V)]^{G}$ which is not preserved by $g \mapsto g^{-t}$.
Let $I \colon \wedge^{4} V^{\vee} \stackrel{\cong}{\to} V$ be the map given by 
\[
(I\omega_{1}, \omega_{2}) = (\omega_{1} \wedge \omega_{2}, \text{vol}_{V}), \ \ \ \ \ \omega_{1} \in \wedge^{4} V^{\vee}, \omega_{2} \in V^{\vee}.
\]
For $\omega_{1}, \ldots, \omega_{5} \in \wedge^{2} V^{\vee}$, define
\[
\Gamma(\omega_{1}, \ldots, \omega_{5}) = (I(\omega_{1} \wedge \omega_{2}) \wedge I(\omega_{3} \wedge \omega_{4}), \omega_{5}).
\]
Extending this map by linearity gives an $\SL(V)$-invariant element $\Gamma \in (\wedge^{2}V)^{\otimes 5}$.
In the basis $\{e_{ij}\}$, we may write $\Gamma$ (up to scale) as follows:
\[
\Gamma = \sum_{\sigma, \sigma' \in S_{5}} (-1)^{|\sigma| + |\sigma'|} e_{\sigma(1)\sigma(2)}\otimes e_{\sigma(3)\sigma(4)} \otimes e_{\sigma'(1)\sigma'(2)} \otimes e_{\sigma'(3)\sigma'(4)} \otimes e_{\sigma(5)\sigma'(5)}.
\]

Let $\phi \colon V \to V^{\vee}$ be the isomorphism induced by the basis $\{e_{i}\}$.
We get an $\SL(V)$-invariant element $\widetilde{\Gamma} := (\wedge^{2}\phi)^{\otimes 5}(\Gamma)$; in the standard basis we have
\[
\widetilde\Gamma = \sum_{\sigma, \sigma' \in S_{5}} (-1)^{|\sigma| + |\sigma'|} e^{\sigma(1)\sigma(2)}\otimes e^{\sigma(3)\sigma(4)} \otimes e^{\sigma'(1)\sigma'(2)} \otimes e^{\sigma'(3)\sigma'(4)} \otimes e^{\sigma(5)\sigma'(5)}.
\]
Define now the function $f \in \CC[\SL(\wedge^{2}(V)]^{G}$ by $f(g) = (\widetilde\Gamma, g\Gamma)$.

Evaluating $f$ at a diagonal matrix $(x_{ij})_{1 < i < j < 5}$, we get
\[
f(x_{ij}) = \sum_{\sigma, \sigma' \in S_{5}} x_{\sigma(1)\sigma(2)} x_{\sigma(3)\sigma(4)}  x_{\sigma'(1)\sigma'(2)}  x_{\sigma'(3)\sigma'(4)}  x_{\sigma(5)\sigma'(5)},
\]
setting $x_{ji} = x_{ij}$ for $j > i$.
This is a homogeneous polynomial of degree 5 which contains $x_{12}^{2}x_{34}x_{35}x_{45}$ as one of its terms.
Such a polynomial cannot be invariant under $x_{ij} \mapsto x_{ij}^{-1}$, even modulo the $\SL(\wedge^{2}V)$-condition $\prod x_{ij} = 1$.
Hence $f$ is not invariant under $g \mapsto g^{-t}$.

\emph{Proof that $g$ is not in the $G$-orbit of $g^{t}$:}
It is enough to construct a function $f \in \CC[\SL(\wedge^{2}V)]^{G}$ which is not preserved by $g \mapsto g^{t}$.
The invariant ring $\CC[\SL(\wedge^{2}V)]^{G}$ equals $\CC[\wedge^{2}V^{\vee} \otimes \wedge^{2}V]^{G}/(\det_{\wedge^{2}V} - 1)$.
We have a decomposition
\[
\CC[\wedge^{2}V^{\vee} \otimes \wedge^{2}V] = \bigoplus_{k} \Sym^{k}(\wedge^{2}V \otimes \wedge^{2}V^{\vee}) = \bigoplus_{k} \bigoplus_{\lambda \vdash k} \SS^{\lambda}(\wedge^{2}V) \otimes \SS^{\lambda}(\wedge^{2}V^{\vee}),
\]
where $\lambda$ denotes a partition and $\SS^{\lambda}$ the associated Schur power \cite[Thm.~2.3.2]{weyman_cohomology_2003}.

The involution $g \mapsto g^{t}$ sends an element $\alpha \otimes \beta \in \SS^{\lambda}(\wedge^{2}V) \otimes \SS^{\lambda}(\wedge^{2}V^{\vee})$ to $\SS^{\lambda}\phi^{-1}(\beta) \otimes \SS^{\lambda}\phi(\alpha)$.
If there is a $\lambda$ such that $\dim \SS^{\lambda}(\wedge^{2}V)^{\SL(V)} > 1$ (and hence $\dim \SS^{\lambda}(\wedge^{2} V^{\vee})^{\SL(V)} > 1$), then a general element $\alpha \otimes \beta \in \SS^{\lambda}(\wedge^{2}V) \otimes \SS^{\lambda}(\wedge^{2}V^{\vee})$ is clearly not preserved by $g \mapsto g^{t}$.

Decomposing $\SS^{\lambda}(\wedge^{2}V)$ as an $\SL(V)$-module is an instance of plethysm.
Running the following lines of Macaulay2 code demonstrates that $\dim \SS^{(5,4,3,2,1)}(\wedge^{2}V)$ contains the $\GL(V)$-representation $(\wedge^{5}V)^{\otimes 6}$ with multiplicity 2, and hence that $\dim \SS^{(5,4,3,2,1)}(\wedge^{2}V)^{SL(V)} = 2$.

\begin{verbatim}
loadPackage "SchurRings";
S = schurRing(QQ,s,5);
plethysm({5,4,3,2,1}, s_{1,1})
\end{verbatim}

So taking general $\alpha \in \SS^{(5,4,3,2,1)}(\wedge^{2}V^{\vee})^{\SL(V)}$ and $\beta \in \SS^{(5,4,3,2,1)}(\wedge^{2}V)^{\SL(V)}$, then $\alpha \otimes \beta \in \CC[\wedge^{2}V^{\vee} \otimes \wedge^{2}V]^{G}$ is not preserved by $g \mapsto g^{t}$.
Let $f \in \CC[\SL(\wedge^{2}V)]^{G}$ be the image of $\alpha \otimes \beta$.
Since $\alpha \otimes \beta$ is homogeneous, it's clear that $f$ is not preserved by $g \mapsto g^{t}$ either.

\emph{Passing from SL to PGL:} In both cases above, letting $f_{\PGL} = f^{2}$, it is easy to see that the function $f_{\PGL}$ on $\SL(\wedge^{2}V)$ descends to $\PGL(\wedge^{2}V)$, is $\PGL(V) \times \PGL(V)^{\op}$-invariant, and is not preserved by $g \mapsto g^{\pm t}$.
\end{proof}

\section{Further questions}
Let $S$ and $Q$ be the universal sub- and quotient bundles on $\Gr(2,5)$, respectively.
As shown in \cite{kapustka,inoue}, the threefold $X_{g}$ admits a smooth degeneration to a variety cut out by a section of $\wedge^{2}Q(1)$.
The family of such degenerate varieties forms a 50-dimensional subset of the 51-dimensional moduli space of deformations of $X_{g}$.
Since $\wedge^{2}Q(1)$ is the normal bundle of $\Gr(2,5)$ in $\PP(\wedge^{2}V)$, these varieties can be thought of as the intersection of $\Gr(2,5)$ with an infinitesimal translate of itself.

Let us recall how these special threefolds come about.
We will work with $\GL(\wedge^{2}V)$ rather than $\PGL(\wedge^{2}V)$.
Consider $g \in \GL(\wedge^{2} V)$ as an element of $\Hom(\wedge^{2}V, \wedge^{2}V)$, and let $s_{g} \in \Gamma(\Gr(2,5), \wedge^{2}V(1)) \cong \wedge^{2}V \otimes \wedge^{2}V^{\vee}$ be the corresponding section.
This gives a section $\wedge^{2} s_{g}$ of $\wedge^{4}V(2)$, and the vanishing locus of $\wedge^{2} s_{g}$ is exactly $X_{g}$.

Consider now $s_{\id} \in \Gamma(\Gr(2,5), \wedge^{2}V(1))$.
We may see $s_{\id}$ as the composition of the canonical and unique section of $\wedge^{2}S(1) \cong \cO_{\Gr(2,5)}$ with the inclusion $\wedge^{2}S(1) \into \wedge^{2} V(1)$.
It follows that $\wedge^{2}s_{\id} = 0$.

Let now $v = g - \id$, so that $s_{v} = s_{g} - s_{\id}$, and consider the 1-parameter family of subvarieties $X_{t}$ of $\Gr(2,5)$ cut out by the sections 
\[
t^{-1}(\wedge^{2}(s_{\id} + ts_{v})) = 2s_{\id} \wedge s_{v} + t\wedge^{2} s_{v}.
\]
We have $X_{1} = X_{g}$, and $X_{0}$ is defined by the section $s_{\id} \wedge s_{v}$ of $\wedge^{4}V (2)$.

We let $Z_{v} = X_{0}$.
Then as shown in \cite{kapustka,inoue}, we may also view $Z_{v}$ as being cut out by the section $p\circ s_{v} \in \Gamma(\Gr(2,5), \wedge^{2}Q(1))$, where $p \colon \wedge^{2}V(1) \to \wedge^{2}Q(1)$ is the obvious map.
The section $p \circ s_{v}$ is transverse, and $Z_{v}$ is smooth for general $v$.

This construction globalises, so that if we let $\cM$ be the blow-up of $\GL(V)$ in the identity, then we have a well-defined family of subvarieties of $\Gr(2,5)$ parameterised by an open subset of $\cM$:
At a point $g \in \cM$ away from the exceptional locus, the corresponding variety is $X_{g}$, and at a point $[v] \in \PP(\gl(\wedge^{2} V))$ in the exceptional locus, it is $Z_{v}$.

The involution of $\GL(\wedge^{2}V)$ given by $g \mapsto g^{-t}$ extends to an involution of $\cM$, which over the exceptional locus $\PP(\gl(\wedge^{2}V))$ is given by $[v] \mapsto [-v^{t}] = [v^{t}]$.
Since the moduli space of Hodge structures is separated, the varieties $Z_{v}$ and $Z_{v^{t}}$ are Hodge equivalent.
The following questions seem natural:
\begin{itemize}
\item Is $Z_{v}$ birational (i.e.~projectively equivalent) to $Z_{v^{t}}$?
\item Is $Z_{v}$ derived equivalent to $Z_{v^{t}}$?
\end{itemize}
A negative answer to the first question would give an alternative proof of Thm.~\ref{thm:XAndYNotBirational}: Indeed, if $X_{g}$ were isomorphic to $X_{g^{-t}}$ for general $g$, that would imply $Z_{v} \cong Z_{v^{t}}$, by Matsusaka--Mumford's theorem \cite{matsusaka_two_1964}.

\bibliographystyle{alpha-abbrv}
\bibliography{bibliography}

\end{document}